\documentclass[reqno, 12pt]{amsart}
\usepackage[pdftex]{graphicx}     

\usepackage{amsfonts}
\usepackage{amssymb,amsmath, color, beton}

\renewcommand{\phi}{\varphi}

\newcommand{\cl}{{\mathcal C}\!\ell }
\newcommand{\ua}{\underline{a}}

\newtheorem{remark}{Remark}
\newtheorem{lemma}{Lemma}
\newtheorem{theorem}{Theorem}
\newtheorem{proposition}{Proposition}
\newtheorem{corollary}{Corollary}
\newtheorem{definition}{Definition}

\begin{document}
\author{Swanhild Bernstein}

\address{Institute of Applied Analysis, University of Mining and Technology, 09599 Freiberg,
Germany, E-mail: swanhild.bernstein@tu-freiberg.de .}

\title[Factorization of  nonlinear equations]{Factorization of the nonlinear Schr\"{o}dinger equation and applications}

\maketitle

\begin{abstract} We consider factorizations of the stationary and non-stationary Schr\"{o}dinger equation in
$\mathbb{R}^n$ which are based on appropriate Dirac operators. These factorizations lead to a Miura transform which
is an analogue of the classical one-dimensional Miura transform but also closely related to the Riccati equation.
In fact, the Miura transform is a nonlinear Dirac equation. We give an iterative procedure which is based on
fix-point principles to solve this nonlinear Dirac equation. The relationship to nonlinear Schr\"{o}dinger equations
like the Gross-Pitaevskii equation are highlighted.\\[2ex]
\end{abstract}

\section{Introduction}
Also today \emph{nonlinear} equations and operators are a challenge. But it is also well-known that transformations
do a good job in understanding and solving nonlinear equations. One famous transform is the Miura transform
(\cite{miura},\,1968) which in one spatial dimension represents a connection between Korteweg--de Vries equations.
The Miura transform can also be considered as a special case of the Riccati equation and can be used to factorize
the (stationary) Schr\"{o}dinger equation. In higher dimensions Clifford algebras and analysis give the possibility to
generalize the factorizations into higher spatial dimensions and even to space-time domains $(t,\,x)\in
\mathbb{R}_+ \times G \subset \mathbb{R}_+ \times \mathbb{R}^n.$\\[1ex]
Factorization of second order differential operators have been studied several times. Factorizations of the
Helmholtz operator in case of a real wave number had been studied by K. G\"{u}rlebeck \cite{g1} and K. G\"{u}rlebeck/W.
Spr\"{o}{\ss}ig in their books \cite{gs1},\,\cite{gs2}. Z.~Xu \cite{xu1},\,\cite{xu2}, F.~Brackx and N.~van~Acker \cite{ba}
and together with R.~Delanghe and F.~Sommen \cite{bads} considered generalized Dirac operators of type $D+k,$ with
the Dirac operator $D$ and $k$ is a complex number. E. Obolashvili \cite{ob1},\,\cite{ob2} and later on Huang Liede
\cite{hl} investigated the case of purly vectorial $k.$ Fundamental solutions for the Dirac operator $D+a,$ where
$a$ is a paravector, are described in \cite{b2}. Factorizations of the Helmholtz equation and the related Dirac
operators $D+\alpha,$ where $\alpha$ is a complex quaternion are investigated by V.~V. Kravchenko and M.~V. Shapiro
\cite{ksh1}, \cite{ksh2}. In the context of Maxwell's equations operators generalized Dirac
operators considered in \cite{McIMi} by A. McIntosh and M. Mitrea.\\[0.5ex]
After this introduction we give a short overview of the Clifford operator calculus in general, starting with
Clifford algebras and then describing Clifford analysis and operator calculus as they needed for applications to
equations of mathematical physics independently of a specific Clifford algebra.\\
In the next section examples are presented, where the main emphasize is laid on the Dirac operator and associated
monogenic functions. The last example are parabolic Dirac operators which allows to factorize the heat equation.\\
Section 3 is concerned with the non-stationary Schr\"{o}dinger operator and appropriate factorizations by parabolic
Dirac operators.\\
In Section 4 we investigate the factorization of the stationary and non-stationary Schr\"{o}dinger equation which will
lead to the Miura transform.\\
With the Miura transform as a nonlinear Dirac equation and a procedure based on fix-point principles to solve this
nonlinear equation is section 5 concerned.\\
The last section deals with the Gross-Pitaevskii equation and it's relation to a nonlinear Schr\"{o}dinger equation and
the Miura transform.

\section{Clifford Operator Calculus}
The connection of analysis and algebra has created an deeper understanding of physics and natural sciences. On way to
algebra to create analysis can be done by the use of Clifford algebras.
\subsection{Clifford algebras} Let $e_1,\,\ldots ,\,e_n$ be an orthonormal basis of
$\mathbb{R}^n.$ The Clifford algebra $\cl _{n,0}$  is the free algebra generated modulo
$$ x^2 = - |x|^2e_0, $$
where $e_0$ is the identity of the Clifford algebra and the basis elements fulfil the anti-commuting relations:
$$ e_ie_j+e_je_i = \pm 2\delta_{ij}e_0, $$
where $\delta_{ij}$ denotes  the Kronecker symbol.
Let $N=\{1,\,2,\,\ldots ,\,n\}$ and for each non-empty subset $A$ of
$N$ set
$$ e_A = e_{A_1}\ldots e_{A_k},\ A=\{A_1,\ldots , A_k\},\ 1\leq A_1<A_2<\ldots <A_k\leq n.$$
By convention $e_{\emptyset} = e_0$ the identity of the Clifford algebra.\\[1ex]
There are three operations on a Clifford algebra. The \emph{principal automorphism}
$$ e_A^{\prime} = (-1)^{|A|}e_A,\quad A\subseteq N, $$
where $|A|$ is the cardinality of $A,$ the \emph{principal anti-automorphism}
$$ e_A^* = (-1)^{\frac{1}{2}|A|(|A|-1)}e_A, $$
while the \emph{conjugation} is the decomposition of the principal automorphism and anti-automorphism
$$ \bar{e}_A = \left(e_A^*\right)^{\prime} = \left(e_A^{\prime}\right)^* = (-1)^{\frac{1}{2}|A|(|A|+1)}e_A. $$
It is easily seen that
$$ \bar{e}_i=-e_i,\,i=1,\,\ldots ,\,n, \quad \mbox{and } \quad \overline{e_ie_j} = \bar{e}_j\bar{e}_i. $$

We will identify the Euclidean space $\mathbb{R}^n$ with $\bigwedge ^1 \cl _{0,n}$ the space of all vectors of $\cl
_{0,n},$ i.e. any element $x\in \mathbb{R}^n$ can be identified with
$$ x = \sum_{i=1}^n x_ie_i. $$
In this way it is easily seen that each non-zero vector has inverse given by $\frac{-x}{|x|^n}.$\\
More on Clifford algebras and Clifford analysis can be found in the pioneering book \cite{bds}, the treatment of
elliptic boundary value problems is considered in \cite{gs1}. The function theoretical basis and the application to
more complicated partial differential equations like the Navier-Stokes equation is contained in \cite{gs2}. A
(comprehensive) function theory for the Dirac operator is the book \cite{dss}. A good description of Clifford
analysis and its connection to classical as well global analysis can be found in \cite{gm}.
\subsubsection{Special cases}\ \\
1. Complex numbers $\mathbb{C}\sim \cl _{0,1}:$ Identifying $1$ with $e_0$ and the complex unit $i$ with $e_1.$ We
obtain that the complex numbers
$\mathbb{C}$ are equivalent to the Clifford algebra $\cl_{0,1}.$  \\[2ex]
2. Quaternions $\mathbb{H}\sim \cl _{0,2}:$ Let us denote the basic elements of the quaternions by $1,\,i,\,j,\,k.$
Then
$$ ij = k = -ji,\,jk = i = -kj,\,ki = j = -ik, $$
and we identify $1=e_0,\,i=e_1,\,j=e_2$ and hence $k=e_1e_2.$ Thus the quaternions $\mathbb{H}$ can be identified with
$\cl_{0,2}.$

\begin{remark} Generally Clifford algebras $\cl _{p,q}$ are defined by the generating elements
$e_i,\,i=1,\,\ldots ,e_n,\,n=p+q,$ where
\begin{eqnarray*}
e_i^2 = -1,\,i=1,\,\ldots ,\,p,\quad e_i^2 = +1,\,i=p+1,\,\ldots ,p+q=n,\\
\quad \mbox{and\ } e_ie_j+e_je_i = -2\delta_{ij}e_0.
\end{eqnarray*}
But we restrict ourself to the cases $\cl_{0,n}$ and $\cl_{n,0}$ because it is possible to factorize all typical types
of second order differential operators by using these Clifford algebras.
\end{remark}

\subsection{Function spaces} \ \\
Let $G\subset \mathbb{R}^n$ be a bounded domain with sufficient smooth boundary $\Gamma = \partial G, $ whose
complement contains a non-empty open set. Then, functions $f$ defined in $G$ with vales in $\cl _{0,n}$ are considered.
These functions may be written as
$$  f(x)= \sum _A f_A\,e_A, \quad x\in G, $$
where $e_A=e_{h_1}\cdots e_{h_k} = e_{h_1\cdots h_k},\ 1\leq h_1 < \ldots h_k \leq n,$ and $e_{\emptyset}=e_0.$
Properties such as continuity, differentiability, integrability and so on, which are ascribed to $f$ have to be
possessed by all components $f_A(x).$ In this way the usual Banach spaces of these functions are denoted by
$C^{\alpha},\ L_p,\ W_p^k.$ The norm in the Banach-spaces $B$ is given by
$$ ||f||_B = \sqrt{\sum_A ||f_A||^2_B }.$$
In particular, we introduce in $L^2(G)$ of \emph{real-valued} function $f_A$ the $\cl _{0,n}$-valued inner product
$$ \langle f,\,g \rangle := \int_{G} \overline{f(x)}\cdot g(x)\,dx, $$
where $\overline{f(x)}= \sum _A f_A\,\overline{e}_A. $ The inner product is related to the $L^2$-norm by
$$ ||f||^2_{L^2} = {\rm sc\,}\langle f,\,f \rangle = \sum_A |f_A|^2(x), $$
where ${\rm sc\,}g = g_0e_0$ denotes the scalar-part of $g= \sum_A g_A.$

\subsection{The operators $D,\,T,\,F$}\ \\
In this section we describe in a general way the basic operators which will become specified by the use of a
specific Clifford algebra. In general, let $G$ be a domain with smooth boundary $\partial G = \Gamma $ and denote
by $n(y)$ the outward pointed normal at $y\in \Gamma.$ There are two ways of consideration, the classical one which
assumes the domain $G$ to have a Ljapunov boundary $\Gamma,$ i.e. $\Gamma$ is of class $C^{1,\alpha}$ with
arbitrary $\alpha >0$ \cite{gs1}, \cite{gs2}. But these strong assumption on the boundary can be weakened to the
consideration of strongly Lipschitz domains \cite{McI1}, \cite{McI2}.

We consider a differential operator $D$ which is a "square root" of a second order differential operator and therefore
possesses a fundamental solution $e(x),$ i.e.
$$De(x)=\delta (x) \quad \mbox{in\ } G.$$
Because the multiplication in Clifford algebras is not commutative we have to distinguish between left-monogenic
functions, i.e. functions $f$ such that $Df=0$ and right-monogenic functions $g,$ which fulfil $gD=0.$  Then the
convolution integral operator
$$Tu= -e\star u = -\int_G e(x-y)u(y)\,dy $$
is a right but not a left inverse of the operator $D$ and there exists an operator
$$(Fu)(x) = \int_{\Gamma} e(x-y)n(y)u(y)\,d\Gamma _y,\ x\in G,$$
such that
$$ DT = I \quad \mbox{and}\quad TD + F = I. $$
The operator $F$ appears by applying the formula of partial integration. By taking limits we obtain
$$\lim_{\ {\mbox{inside\,}G/{\mbox{outside\,}G \to \Gamma }}} Fg  = \pm g  + Sg, $$
with a Hilbert-type operator
$$(Su)(x) = 2\int_{\Gamma}e(x-y)n(y)u(y)\,d\Gamma _y, x\in \Gamma .$$

\subsection{Decompositions}
Another general principal is that
\begin{eqnarray} \label{gl1}
 L^2(G)= {\rm ker \,}D \oplus \overline{{\rm im \,} D^*.}
\end{eqnarray}
It is a challenge to prove that the (orthogonal) complement is ${\rm im \,} D^*,$ this set is closed itself. Now, it is
easy to see that this decomposition defines (orthogonal) projections $\mathbb{P}$ and $\mathbb{Q}.$
\begin{remark} The just mentioned decomposition are orthogonal in $L^2$ and can be extended to direct sums in
$L^p,\,1<p<\infty,$ and in Sobolev spaces $W_p^k,\,p\geq 1.$
\end{remark}
The connection between the operator $F$ and the Hilbert-type operator $S$ give raise to another decomposition which is
based on the (orthogonal) projections $P=\frac{1}{2}(I+S)$ and $Q=\frac{1}{2}(I-S).$  This is a Hardy-type
decomposition
\begin{eqnarray} \label{gl2}
 L^2(\Gamma )= {\rm im \,} P \oplus {\rm im \,} Q .
\end{eqnarray}
An interesting connection is given between (\ref{gl1}) and (\ref{gl2}) in the way that $f\in {\rm ker \,}D$ if and only
if ${\rm tr\,}f \in {\rm im \,} P.$ This connection is given by the Plemelj-Sokhitzkij's formulas which tell us that
${\rm tr\,}f \in {\rm im \,} P$ if and only ${\rm tr\,}f$ has a monogenic continuation $f$ into $G,$ i.e. $Df=0.$

\subsection{Examples}

\subsubsection{Dirac operator}
Here, $D$ denotes the Dirac operator
$$ Df=\sum_{j=1}^{n} e_j\frac{\partial f}{\partial x_j} , $$ which factorizes the Laplacian:
$$ DD = -\Delta_n, $$
where $\Delta_n $ denotes the Laplacian in $\mathbb{R}^n,$ and we obtain the fundamental solution of the Dirac operator
by applying the Dirac operator to the fundamental solution of the Laplacian $H.$ We get the generalized Cauchy kernel
$$ e(x)=\frac{-1}{\omega_n}\frac{x}{|x|^{n}}, $$
where $\omega_n$ denotes the surface area of the unit sphere in $\mathbb{R}^{n}.$ The first conclusion is
\begin{corollary}[Borel-Pompeiu formula]\label{bp}
Let $G$ be a bounded strongly Lipschitz domain of $\mathbb{R}^n.$ If $f\in W^1_p(G),\,1<p<\infty,$ we have
$$ Ff + TDf = \left\{ \begin{array}{cl} f , & \mbox{in\,} G ,\\ 0, & \mbox{in\,} \mathbb{R}^n\backslash \bar{G}.
\end{array} \right. $$
\end{corollary}
Similar to the case of one-dimensional complex variables we can define Hardy spaces \cite{mitrea1}. Let us denote
by $G_+=G$ and $G_-=\mathbb{R}^n\backslash \bar{G}.$
\begin{definition}
For a function $F:\,G_{\pm} \to \cl_{0,n}$ and $1<p<\infty $ we set
$$ ||F||_{\mathcal{H}^p}:=\sup _{\delta >0} \left\{ \int_{\Gamma} |F(x\pm \delta)|^p\,d\Gamma_x \right\}^{\frac{1}{p}},
$$
and define the \emph{Hardy spaces} of monogenic functions
$$ \mathcal{H}^p(G_{\pm}):=\{ f \mbox{\ is left monogenic in\ } G_{\pm}; ||f||_{\mathcal{H}^p}< \infty \}. $$
\end{definition}
There is also a connection to the boundary values of functions from the Hardy spaces $\mathcal{H}^p(G_{\pm}).$

\begin{corollary} Let $1<p<\infty. $ For a left monogenic function $H$ in $G,$ the following are equivalent
\begin{enumerate}
    \item $ H \in \mathcal{H}^p(G);$
    \item There exists $h\in L^p(\Gamma )$ such that $H$ is the (right) Cauchy integral extension of $h,$ i.e.
    $ H(x)= Fh(x) $ for $x$ in $G;$
    \item $H$ has a non-tangential boundary limit $H^+(x)$ at almost any point $x\in \Gamma ,$ i.e. there exists
    $$ H^+(x):=\lim_{\{y\in x+\Sigma _{\alpha},\,y\to x\}} H(y) \mbox{\ for a.e.\ } x\in \Gamma, $$
    and $H$ is the (right) Cauchy integral extension of its boundary trace.
    In addition the norms $||\cdot
    ||_{\mathcal{H}^p}$ and $|| H^+||_{L^p}$ are equivalent, i.e.
    $$ ||H||_{\mathcal{H}^p} = || H^+||_{L^p}.$$
\end{enumerate}
\end{corollary}
The most important singular convolution operator is the singular Cauchy integral operator $S$ on $\Gamma ,$ defined for
almost all $x\in \Gamma $ by
$$ Sf(x)=2 \lim_{\varepsilon \to 0} \int_{\{y\in\Gamma :\,|x-y|> \varepsilon \}} e(x-y)n(y)f(y)\,d\Gamma _y. $$
\begin{corollary}[Plemelj-Sokhotzkij's formula]
Let $G$ be a bounded strongly Lipschitz domain of $\mathbb{R}^n$ with boundary
$\partial G = \Gamma.$ If $f\in L^p(\Gamma ),\,1<p<\infty,$ then
\begin{align*}
(Pf)(x)= & \lim_{\{\varepsilon \to 0} \int_{y\in\Gamma:\,|x-y|>\varepsilon\}} e(x-y)n(y)f(y)\,d\Gamma _y +
\frac{1}{2}f(x)
\\ = & \frac{1}{2}\left(I+S\right)f(x).
\end{align*}
and hence
$$ (Qf)(x)=\frac{1}{2}\left(I-S\right)f(x). $$
\end{corollary}
Moreover, we have
$$ Sf = 2Pf(x)- f(x) \quad \mbox{and} \quad  Sf = -2Qf + f. $$
The equality of the integrals is shown in \cite{McI1}. It follows from Calder\'{o}n-Zygmund theory. We have
\begin{corollary} The space $L^p(\Gamma ),\, 1<p<\infty,$ posses the  decomposition
$$ L^p(\Gamma ) = {\rm im\,} P\cap L^p(\Gamma ) \oplus {\rm im\,} Q \cap L^p(\Gamma ), $$
where the decomposition is orthogonal in case $p=2$ and direct in case $p\not = 2.$ Moreover, the spaces ${\rm im\,}
P\cap L^p(\Gamma )$ and ${\rm im\,} Q \cap L^p(\Gamma )$ can be identified with the boundary values of the Hardy spaces
$\mathcal{H}^p(G_{\pm}).$
\end{corollary}
The Teodorescu transform is given by
$$ Tf(x)= \frac{1}{\omega_n}\int_G\frac{x-y}{|x-y|^{n}}f(y)\,dy $$
and from Borel-Pompeiu's formula we immediately conclude that
$$ {\rm tr\,} Tf \in {\rm im\,}Q. $$

\subsubsection{Cauchy-Riemann operator}
Now, $D$ denotes the Cauchy-Riemann operator
$$ Df=\frac{\partial f}{\partial x_0} + \sum_{j=1}^{n} e_j\frac{\partial f}{\partial x_j} \quad \mbox{and} \quad
\bar{D}f=\frac{\partial f}{\partial x_0} - \sum_{j=1}^{n} e_j\frac{\partial f}{\partial x_j} . $$ Also the
Cauchy-Riemann operator factorizes the Laplacian. We have
$$ D\overline{D} = \Delta _{n+1}, $$
where $\Delta_{n+1} $ denotes the Laplacian in $\mathbb{R}^{n+1},$ and we obtain the fundamental solution of the
generalized Cauchy kernel
$$ e(x)=\frac{1}{\omega_{n+1}}\frac{\overline{x}}{|x|^{n+1}}, $$
where $\omega_{n+1}$ denotes the surface area of the unit sphere in $\mathbb{R}^{n+1}.$ \\
The Cauchy-Riemann operator can be transform into a Dirac operator by multiplication with $e_{n+1}:$
$$ \frac{\partial f}{\partial x_0} - \sum_{j=1}^{n} e_j\frac{\partial f}{\partial x_j} e_{n+1} =
   \frac{\partial f}{\partial x_0}e_{n+1} - \sum_{j=1}^{n} e_je_{n+1}\frac{\partial f}{\partial x_j} =
   \frac{\partial f}{\partial x_0}E_{n+1} - \sum_{j=1}^{n} E_j\frac{\partial f}{\partial x_j}. $$
On the other hand the Dirac operator can be transformed into the Cauchy-Riemann operator:
$$ \sum_{j=1}^{n} e_j\frac{\partial f}{\partial x_j} (-e_{n}) =
   \frac{\partial f}{\partial x_n} + \sum_{j=1}^{n-1} ( -e_je_{n})\frac{\partial f}{\partial x_j} =
   \frac{\partial f}{\partial x_n} + \sum_{j=1}^{n-1} E_j\frac{\partial f}{\partial x_j}. $$

\subsubsection{$\overline{\partial}$-operator} The simplest case occurs when the underlying Clifford algebra is just
equivalent to the complex numbers, i.e. $\cl_{0,1}\sim \mathbb{C}.$ In this case we have the complex Cauchy-Riemann
operators
$$ \bar{\partial} _z = \frac{\partial}{\partial x} + i \frac{\partial}{\partial
y} \quad \mbox{and} \quad \partial _z = \frac{\partial}{\partial x} - i \frac{\partial}{\partial y} $$ which can be
identified with
$$ D  = \frac{\partial}{\partial x} + e_1\frac{\partial}{\partial
y}\quad \mbox{and}\quad  \bar{D} = \frac{\partial}{\partial x} - e_1 \frac{\partial}{\partial y}
$$
where $\bar{\partial}$ is identified with $D$ and we have the properties
$$ \bar{\partial}\partial = \partial\bar{\partial} = D\bar{D} = \bar{D}D = \Delta = \frac{\partial^2}{\partial x^2} +
\frac{\partial^2}{\partial y^2}.$$
The fundamental solution of the $\bar{\partial}$-operator is the Cauchy kernel
$$ e(z)=-\frac{1}{2\pi\,i}\frac{1}{z} = \frac{1}{2\pi}\frac{1}{z}i \sim \frac{1}{2\pi}\frac{\overline{z}}{|z|^2}e_1
$$
and we have the generalized Cauchy formula (Cauchy-Green formula, Borel formula, ...) as an anolog of the
Borel-Pompeiu formula:
$$ f(z)=\frac{1}{2\pi\,i}\int_{\Gamma} \frac{f(\zeta)}{\zeta - z}\,d\zeta - \frac{1}{\pi}\int_G
\frac{\bar{\partial}f(\zeta)}{\zeta -z},\,d\zeta .$$

\subsubsection{Generalized Dirac and Cauchy-Riemann operators}
We start with the Dirac operator and consider the disturbed operator
$$ D + \alpha = \sum_{j=1}^{n} e_j\frac{\partial f}{\partial x_j} + \alpha $$ and
$$ (D + \alpha)(D - \alpha)= DD - \alpha ^2 = -\Delta -\alpha ^2, $$ where $\alpha$ is a scalar-valued function. But we
can do even more
$$ (D + a)(D+\bar{a})= DD+ a\bar{a} = -\Delta - |a|^2. $$
An extensive study of these generalized or disturbed Dirac operators can be found in \cite{ksh1} and \cite{ksh2}.
\subsubsection{Parabolic Dirac operators}
Recently (cf. \cite{cks}), it was also possible to factorize the heat equation by a Witt-type basis. Which means to
add two additional generating elements $\frak{f}$ and $\frak{f}^+,$ where
$$ \begin{array}{rcl} \frak{f}\frak{f}^+ +\frak{f}^+ \frak{f} & = & 1, \\[0.5ex]
                          \frak{f}^2 = \left(\frak{f}^+\right)^2 & = & 0, \\[0.5ex]
                          \frak{f}\,e_j + e_j\,\frak{f} & = & 0, \\[0.5ex]
                           \frak{f}^+ e_j + e_j\,\frak{f}^+ & = & 0. \end{array} $$
and
$$ {\bf\sl D}_{x,t}^{\pm} = \sum_{j=1}^{n} e_j\frac{\partial f}{\partial x_j} + \frak{f}\partial _t \pm \frak{f}^+, $$
where $$ ({\bf\sl D}_{x,t}^{\pm})^2 = -\Delta \pm \partial _t, $$ i.e. ${\bf\sl D}_{x,t}^+$ factorizes the heat operator.
The fundamental solution for the heat operator is
$$ e(x,t)=\frac{H(t)}{(2\sqrt{\pi\,t})^n}\,e^{-\frac{|x|^2}{4t}}, $$
where $H(t)$ denotes the Heaviside-function. More on parabolic Dirac operators and the application to the
non-stationary Navier-Stokes equation can be found in the already mentioned paper \cite{cks}.

\section{Parabolic Dirac operators related to the Schr\"{o}dinger equation}
Let us denote with $C$ the space-time domain $C=G_t\times [0,\,\infty ) \subset \{(x,t)\in \mathbb{R}^n\times
\mathbb{R}^+\}$ and by $\partial C = G\cup M,$ where $M=\Gamma _t \times [0,\infty),$  its sufficiently smooth
boundary. With the subindex $t$ we emphasis the fact that the space domain (and therefore $M$) can depend on the time
level $t,$ but has to be bounded for all times $t\in [0,\,\infty ).$

The consideration of the Schr\"{o}dinger equation is similar to the heat equation but the Schr\"{o}dinger equations reads as
$$ -i\partial _t - \Delta_n , $$
where $\Delta _n $ again denotes the Laplacian in $\mathbb{R}^n.$ It is easily seen that we obtain a fundamental
solution of the Schr\"{o}dinger equation from the fundamental solution of the heat equation by the transform $t\to it, $
hence
$$ E(x,t)=\frac{H(t)}{(2\sqrt{\pi\,i\,t})^n}\,e^{i\frac{|x|^2}{4t}}, $$
where $H(t)$ denotes the Heaviside-function.

We consider the operators
$$ D_{x,t}^{\pm} = \sum_{j=1}^{n} e_j\frac{\partial f}{\partial x_j} +\frak{f}\,\partial _t \pm i\frak{f}^+,$$
and hence
\begin{multline*}
(D_{x,t}^{\pm})^2u = D_{x,t}^{\pm}D_{x,t}^{\pm}u =(D_x+\frak{f}\partial_t \pm i\frak{f}^+)
(D_xu+\frak{f}\partial_tu \pm i\frak{f}^+u)\\
=D_xD_xu - \frak{f}D_x\partial_tu \mp i\frak{f}^+D_xu + \frak{f}D_x\partial_tu \\ + \frak{f}^2\partial_t^2u \pm
\frak{f}\frak{f}^+\partial_tu
\pm i\frak{f}^+D_xu +\pm i\frak{f}^+\frak{f}\partial_tu - \left(\frak{f}^+\right)^2u \\
= -\Delta u \pm i(\frak{f}\frak{f}^+ + i\frak{f}^+\frak{f})\partial _t u = -\Delta u \pm i\partial_t u.
\end{multline*}

Next, we consider the reduced differential operator
$$ D_{x,t} = \sum_{j=1}^{n} e_j\frac{\partial f}{\partial x_j} + \frak{f}\,\partial _t .$$
If we introduce the sigma-form
$$ d\sigma _{x,t} = D_{x,t}\rfloor  dV_x\,dt $$
we write Stokes theorem as follows
$$ \int_{\partial C} f\,d\sigma _{x,t}\,g = \int_C (f\,D_{x,t})\,g + f\,(D_{x,t}\,g) \, dV_x\,dt $$
and thus
$$ \int_{\partial C} f\,d\sigma _{x,t}\,g = \int_C (f\,D^+_{x,t})\,g + f\,(D^-_{x,t}\,g) \, dV_x\,dt. $$
From the fundamental solution of the Schr\"{o}dinger equation we can switch to the fundamental solution of $-\Delta +
i\partial _t$ by making a reflection $t\to -t$ and applying $D_{x,t}^+$ from the right we obtain the fundamental
solution
\begin{multline*}
e(x,t)=E(x,-t)D_{x,t}^+  =
\frac{H(-t)}{(2\sqrt{\pi\,i\,(-t)})^n}\,e^{-i\frac{|x|^2}{4t}}(D_x+\frak{f}\partial _t +i \frak{f}^+ ) \\
= \frac{H(-t)}{(2\sqrt{\pi\,i\,(-t)})^n}\,e^{-i\frac{|x|^2}{4t}}\left(\frac{-i}{2t}\sum_{j=1}^nx_je_j +
\frak{f}\left(-\frac{n}{2t} + \frac{i|x|^2}{4t^2}\right)+i\frak{f}^+\right).
\end{multline*}
which gives the \emph{Borel-Pompeiu formula} for the operator $D_{x,t}^-:$
$$ \int_{\partial C} e(x-x_0,\,t-t_0)\,d\sigma_{x,t}g(x,t) = g(x_0,t_0) + \int_C e(x-x_0,\,t-t_0)(D^-_{x,t}\,g) \,
dV_x\,dt $$ and for $g\in {\rm ker\,}D_{x,t}^- $ Cauchy's integral formula:
$$ \int_{\partial C} e(x-x_0,\,t-t_0)\,d\sigma_{x,t}g(x,t) = g(x_0,t_0) .$$

\section{Factorization of the Schr\"{o}dinger equation}
The factorization of the Schr\"{o}dinger equation is treated in several papers. Mostly the factorization is related to
systems of differential equations in mathematical physics \cite{b1}, \cite{kk1}, \cite{kk2}. In \cite{kkw}
especially the relationship to the Riccati equation is investigated and quaternionic generalizations of the Riccati
equation had are established.\\
A detailed factorization of the stationary Schr\"{o}dinger equation for $n=3,4$ in an quaternionic context was treated
in \cite{bg}.\\
The Miura transform can be obtained by factorizing the 1D Schr\"{o}dinger equation:
$$ -\frac{d^2}{dx^2}u-v(x)u = -\left(\frac{d}{dx}+a (x)\right)\left(\frac{d}{dx}-a (x)\right) u. $$
We have
\begin{eqnarray*}
-\left(\frac{d}{dx}+a (x)\right)\left(\frac{d}{dx}-a (x)\right) u \\
= - \left(\frac{d^2u}{dx^2}-\frac{d}{dx}(a(x)u)+ a(x)\frac{du}{dx} - a^2(x)u \right)  \\
= -\left(\frac{d^2u}{dx^2}+\left(-\frac{da}{dx}-a^2(x)\right)u\right)
\end{eqnarray*}
and thus
$$ \frac{da}{dx}+a^2(x)=-v(x).$$

\subsection{Factorization of the stationary Schr\"{o}dinger equation}
The Helmholtz equation can be treated quite similar to the Laplace equation and can be used to describe Maxwell's
equations. The relationship between Maxwell's equations as a system of first order differential equations and the
formulation in terms of the Helmholtz equation is nicely characterized by the factorization of the Helmholtz
equation. The Helmholtz equation and the factorization in generalized Dirac operators is discussed in \cite{ksh1},
\cite{ksh2}, \cite{kg} and \cite{mitrea2}. In higher dimensions it is quite clear that the Helmholtz operator may
be factorized by using Dirac operators:
$$ -\Delta - k^2 = (D+k)(D-k). $$
More difficult is the case of a variable potential $V_0(x).$ We consider the case $u=u(x)=u(x)e_0$ and look for
suitable functions $a $ with
\begin{align*}
(-\Delta -V_0(x))u & = (D+a (x))(D-a (x))u \\
                   & = DDu-D(a (x)u)+a (x)Du -a ^2(x)u \\
                   & = -\Delta u  \underline{-Du\cdot a(x)+a(x)Du }-Da(x)\cdot u - a ^2(x)u.
\end{align*}
The underlined part does not vanish, because of the non-commutativity of the multiplications of elements of a Clifford
algebra. Thus, we will change our approach using a multiplication operator $M^{a (x)}$ defined by
$$ M^{a (x)}u(x):= u(x)\cdot a (x). $$
Therefore
\begin{align*}
(-\Delta -V_0(x))u & = \left(D+M^{a (x)}\right)\left(D-M^{a (x)}\right)u \\
                   & = DDu - D(ua(x)) + Du\cdot a(x) - uDa(x) -u a^2(x) \\
                   & = -\Delta u -\left(Da(x) +a^2(x) \right)u
\end{align*}
or
\begin{eqnarray} \label{m}
 Da(x) +a^2(x) = V_0(x).
\end{eqnarray}
Equation (\ref{m}) will be called \emph{Miura transform}. It is nonlinear Clifford-valued first-order partial
differential equation or equivalently a nonlinear real-valued first-order \emph{system} of partial differential
equations.

\subsection{Factorization of the non-stationary Schr\"{o}dinger equation}
Now, by using the operator $D_{x,t}^-$ we can also obtain a factorization of the \emph{non-stationary} Schr\"{o}dinger
equation. Let us start with
\begin{multline*} \left(D_{x,t}^-+M^{a (x,t)}\right)\left(D_{x,t}^--M^{a (x,t)}\right)u(x,t) \\
=D_{x,t}^-D_{x,t}^-u(x,t)-D_{x,t}^-(u(x,t)a(x,t))+(D_{x,t}^-u(x,t))a(x,t) -u(x,t)a^2(x,t).
\end{multline*}
We have to evaluate
\begin{eqnarray*}
D_{x,t}^-(u(x,t)a(x,t)) = \left(D_x+\frak{f}\partial _t - i\frak{f}^+\right)(u(x,t)a(x,t)) \\
= D_x(u(x,t)a(x,t)) + \frak{f}(\partial _t(u(x,t)a(x,t))) -i \frak{f}^+(u(x,t)a(x,t))
\end{eqnarray*}
If we again assume that
$u=u(x,t)=u(x,t)e_0$ is a scalar-valued function we can conclude
\begin{eqnarray*}
= (D_{x}u(x,t))a(x,t) +
u(x,t)(D_{x}a(x,t))+\frak{f}(\partial _tu)a(x,t) +u(x,t)\frak{f}(\partial _ta(x,t))\\
-i\frak{f}^+(u(x,t) a(x,t))
\end{eqnarray*}
and finally
\begin{multline*}
\left(D_{x,t}^-+M^{a (x,t)}\right)\left(D_{x,t}^--M^{a (x,t)}\right)u(x,t) =\\
=D_{x,t}^-D_{x,t}^-u(x,t) -(D_{x}u(x,t))a(x,t) -
u(x,t)(D_{x}a(x,t))-\frak{f}(\partial _tu(x,t))a(x,t) \\
- u(x,t)\frak{f} (\partial _ta(x,t))+
i\frak{f}^+(u(x,t)a(x,t))+(D_{x,t}^-u(x,t))a(x,t) -u(x,t)a^2(x,t) \\
=(-i\partial_t-\Delta_n)u(x,t)-u(x,t)\left(D_{x}a(x,t)+\frak{f} (\partial _ta(x,t))+a^2(x,t)\right)
\end{multline*}
If we compare the last relation with the non-stationary Schr\"{o}dinger
equation we obtain
\begin{eqnarray} \label{eqpot}
D_{x}a(x,t)+\frak{f}(\partial_ta(x,t))+a^2(x,t)=V(x,t)
\end{eqnarray}
This equation is similar to the equation
\begin{eqnarray} \label{kk}
 \partial _tg + Dg + |g|^2 = 0
\end{eqnarray}
which has formal similarities to the canonical Riccati equation. Equation (\ref{kk} was investigated in \cite{kk3}.
We specify the equation (\ref{eqpot}) in the following way. In applications the potential usually depends only on
$x$ and not on $t.$ Therefore, if we consider a potential $V=V(x)$ equation (\ref{eqpot}) suggest to consider also
functions $a =a(x)$ which depend only on $x.$ Which implies that we will consider vector-valued functions $\ua (x)$
depending only on $x$ and we end up with the equation
\begin{eqnarray} \label{eqnl}
D_{x}\ua(x)+\ua^2(x)=V(x),
\end{eqnarray}
which is the same as for the factorization of the stationary Schr\"{o}dinger equation.
\section{The nonlinear Dirac equation}
In this section we want to solve equation (\ref{eqnl}) under the following assumptions:
\begin{enumerate}
    \item the potential $V=V(x)$ is a scalar-valued function that depends only on $x\in G,$  where $G$ is a
    compact and smooth domain of $\mathbb{R}^n.$
    \item we are looking for a vector-valued function $\ua(x)=\sum_{i=1}^n a_i(x)e_i $ that depends only on
    $x\in G .$
\end{enumerate}

If we apply the $T$-operator and use Borel-Pompeiu's formula (cf. Corollary \ref{bp}) we obtain
$$ \ua  - F\ua = T(V+|\ua|^2), $$
because we are only looking for $\ua \in {\rm im\,}Q$ we have $F\ua = 0$ and we end up with
$$ \ua = T(V+|\ua|^2) $$
Our iteration procedure reads now as follows:
$$ \ua_n:=T(V+|\ua_{n-1}|^2), \quad n=1,2,\ldots , $$
Obviously, $V(x)+|\ua_{n-1}|^2\in \mathbb{C}$ reproduces $\ua_n=T(V+|\ua_{n-1}|^2)\in \bigwedge ^1 \cl _{0,n}$ and
because of ${\rm tr\,}T \in {\rm im\,}Q,$ also ${\rm tr\,}\ua_{n}\in {\rm im\,}Q.$\\[0.5ex]
That means that we also have to assume ${\rm tr\,}\ua_{n-1}\in {\rm im\,}Q.$\\[0.5ex]
Now, we consider the regularity. The Teodorescu transform is a weakly singular integral operator and therefore
$$ V \in L^p(G) \Rightarrow TV\in W^1_p(G),\ 1<p<\infty, $$
on the other hand Sobolev's embedding theorems $ W^1_p(G)  \hookrightarrow L^{p^*}(G),$ where $p^*\leq
\frac{np}{n-p}\leq n,$ with embedding constant $C$ leads to the following. Set $p^*=2p$ and assume that $
\ua_{n-1}\in W^1_p(G).$ Then
\begin{align*}
||\,|\ua_{n-1}|^2\,||_{L^{p}} & = \left( \int_G |\ua_{n-1}(x)|^{2p} \,dG \right)^{1/p}
                               = \left( \int_G \left(\sum_{j=1}^n a_{j,n-1}^2(x)\right)^{\frac{1}{2}\cdot 2p}\,dG \right)^{1/p}\\
& = \left(\int_G \left(\sum_{j=1}^n a_{j,n-1}^2(x)\right)^{p}\,dG \right)^{1/p} = ||\sum_{j=1}^n a_{j,n-1}^2(x)
||_{L^p} \\
     & \leq \sum_{j=1}^n ||a^2_{j,n-1}(x)||_{L^p} = \sum_{j=1}^n \left(\int_G a_{j,n-1}^{2p}(x)\,dG\right)^{1/p}\\
     & = \sum_{j=1}^n ||a_{j,n-1}(x)||_{L^{2p}}^2 \leq  \sum_{j=1}^n  C^2\,||a_{j,n-1}(x)||_{W_p^1}^2 =
     C^2\,||\ua_{n-1}(x)||_{W_p^1}^2,
\end{align*}
i.e. $\ua_{n-1}\in W^1_p(G)$ implies $|\ua_{n-1}|^2\in L_p(G)$ and finally also
$T|\ua_{n-1}|^2\in W^1_p(G).$\\[0.5ex]
Thus we assume that $n>p$ and
$$ \frac{np}{n-p} \geq 2p \iff np\geq 2np-p^2 \iff 2p^2 \geq np \iff 2p\geq n \iff p\geq \frac{n}{2}. $$
Therefore we need the following assumptions too:
\begin{enumerate}
    \item $n>p\geq \frac{n}{2}>1,$ if $n>2$ (or $n>p>1$ if $n=2$), $p\in \mathbb{R},$
    \item $\ua_{n-1} \in W_p^1(G)$ and ${\rm tr\,}\ua_{n-1} \in {\rm im\,}Q\cap W_p^{1-\frac{1}{p}}(\Gamma ).$
\end{enumerate}
The last assumption involves that if $\ua_{n-1} $ belongs to $W_p^1(G)$ due to the trace theorem for Sobolev spaces
${\rm tr\,}\ua _{n-1}$ belongs to $W_p^{1-\frac{1}{p}}(\Gamma ).$ \\[0.5ex]
In order to apply  fixed-point theorems we first prove the boundedness of the sequence $\{\ua_n\}.$ We have
\begin{align}
||\ua_n||_{W_p^1} & \leq ||T||_{[L^p,\,W_p^1]}\left(||V||_{L^p}+||\,|\ua_{n-1}|^2||_{L^p}\right) \nonumber \\
                  & \leq k_1\left(||V||_{L^p}+C^2||\ua_{n-1}||^2_{W_p^1}\right) \nonumber \\
                  & \leq k_1||V||_{L^p}+k_2||\ua_{n-1}||^2_{W_p^1}, \label{f1}
\end{align}
where $k_1=||T||_{[L^p,\,W_p^1]},\ k_2=k_1C^2,$ and $C$ is again the embedding constant for $ W^1_p  \hookrightarrow
L^{2p}.$
\begin{lemma} \label{lem1}
We assume that $ ||V||_{L^p}\leq \frac{1}{4k_1k_2}$ and denote by
\begin{eqnarray} \label{f2}
W=\sqrt{\frac{1}{4k_2^2}-\frac{k_1}{k_2}||V||_{L^p}}.
\end{eqnarray}
Then
$$ \frac{1}{2k_2}-W\leq ||\ua_{n_0}||_{W_p^1} \leq \frac{1}{2k_2}+W $$
implies that $||\ua_{n_0+1}||_{W_p^1}\leq ||\ua_{n_0}||_{W_p^1}$ and the sequence $\{\ua_n\}_{n\geq n_0}$ will be
bounded by $||\ua_{n_0}||_{W_p^1}.$
\end{lemma}
\begin{proof} The inequality (\ref{f2}) can be rewritten as
$$ -W \leq ||\ua_{n_0}||_{W_p^1} - \frac{1}{2k_2} \leq W $$
and hence is equivalent to
$$ ||\ua_{n_0}||^2_{W_p^1}-\frac{1}{k_2}||\ua_{n_0}||_{W_p^1}+\frac{1}{4k_2^2} \leq W^2 =
\frac{1}{4k_2^2}-\frac{k_1}{k_2}||V||_{L^p}$$ and we obtain that
$$ ||\ua_{n_0}||^2_{W_p^1}-\frac{1}{k_2}||\ua_{n_0}||_{W_p^1} + \frac{k_1}{k_2}||V||_{L^p} \leq 0 $$
and using (\ref{f1})
$$ ||\ua_{n_0+1}||_{W^1_p}\leq k_1||V||_{L^p}+ k_2 ||\ua_{n_0}||^2_{W_p^1}   \leq ||\ua_{n_0}||_{W_p^1}.$$
\end{proof}
\begin{theorem}
Suppose that $n>p\geq \frac{n}{2}>1,\ p\in \mathbb{R},$ and $||V||_{L^p}\leq \frac{1}{4k_1k_2}.$ Then the equation
$$ \ua = T(V+|\ua|^2) $$
has at least one solution with $$ ||\ua||_{W^1_p} \leq \frac{1}{2k_2}+
\sqrt{\frac{1}{4k_2^2}-\frac{k_1}{k_2}||V||_{L^p}} . $$
\end{theorem}
\begin{proof} From \ref{lem1} we immediately see that
$$  ||\ua_{n-1}||_{W_p^1} \leq \frac{1}{2k_2}+W \Rightarrow ||\ua_{n}||_{W_p^1} \leq \frac{1}{2k_2}+W .$$
If we start with $\ua_0 \in W_p^1(G)$ such that $||\ua_{0}||_{W_p^1} \leq \frac{1}{2k_2}+W $ then the sequence
$\{\ua_n\}_{n\in\mathbb{N}}$ is bounded from above by $\frac{1}{2k_2}+W.$ Hence there exists a subsequence
$\{\ua_{n^{\prime}}\} \subset W^1_p(G)$ with $\ua_{n^{\prime}}\rightharpoonup \ua$ as $n^{\prime}\to \infty .$ Due to
the continuity of $T:\,L^p(G)\to W_p^1(G)$ we conclude
$$ \ua = T(V+|\ua|^2). $$
The norm estimates follows from the weak convergence of $\ua_{n^{\prime}}$ in a convex set.
\end{proof}
But we can do better
\begin{theorem}
We assume that
\begin{enumerate}
    \item $n>p\geq \frac{n}{2}>1,\ p\in \mathbb{R},$
    \item $||V||_{L^p}\leq \frac{1}{4k_1k_2},$
    \item $\ua_0\in W_p^1(G)$ and ${\rm tr\,}\ua_0\in {\rm im\,}Q\cap W_p^{1-1/p}(\Gamma),$
    \item $||\ua||_{W^1_p} \leq \frac{1}{2k_2}-
\sqrt{\frac{1}{4k_2^2}-\frac{k_1}{k_2}||V||_{L^p}}.$
\end{enumerate}
and define the sequence $\{\ua_n\}_{n\in\mathbb{N}} $ by
$$ \ua_n = T(V+|\ua_{n-1}|^2),\quad n=1,2,\ldots .$$
Then there exists a \emph{unique} solution $\ua \in W_p^1(G)$ with ${\rm tr\,}\ua \in {\rm im\,}Q\cap
W_p^{1-1/p}(\Gamma)$ of
$$ \ua = T(V+|\ua|^2) .$$
Moreover, the sequence $\{\ua _n\}_{n\in\mathbb{N}} $ converges to $\ua$ in $W_p^1$ and the solution $\ua$ fulfils the
estimate
$$ ||\ua||_{W_p^1} \leq \frac{1}{2k_2}- \sqrt{\frac{1}{4k_2^2}-\frac{k_1}{k_2}||V||_{L^p}}. $$
\end{theorem}
\begin{proof} To apply Banach's fix-point theorem we investigate the contractivity of the mapping $T(v+|\ua|^2).$
At first we get
$$ ||\ua_n-\ua_{n-1}||_{W_p^1} = ||T(|\ua_{n-1}|^2-|\ua_{n-2}|^2)||_{W_p^1}
\leq k_1\,||\,|\ua_{n-1}|^2-|\ua_{n-2}|^2\,||_{L^p}.$$
Furthermore,
\begin{align*}
||\,|\ua_{n-1}|^2-|\ua_{n-2}|^2\,||_{L^p} & = ||(|\ua_{n-1}|-|\ua_{n-2}|)(|\ua_{n-1}|+|\ua_{n-2}|)||_{L^p} \\
     & \leq ||\,|\ua_{n-1}|-|\ua_{n-2}|\,||_{L^{2p}}||\,|\ua_{n-1}|+|\ua_{n-2}|\,||_{L^{2p}} \\
     & \leq ||\,|\ua_{n-1}|-|\ua_{n-2}|\,||_{L^{2p}}||\,|\ua_{n-1}|^2||^{1/2}_{L^{p}}+||\,|\ua_{n-2}|^2||^{1/2}_{L^{p}} \\
& \leq ||\,|\ua_{n-1}|-|\ua_{n-2}|\,||_{L^{2p}}\left(C||\ua_{n-1}||_{W^1_p}+
      C||\ua_{n-2}||_{W^1_p}\right) \\
     & \leq 2C||\,|\ua_{n-1}|-|\ua_{n-2}|\,||_{L^{2p}}\left(\frac{1}{2k_2}-
     \sqrt{\frac{1}{4k_2^2}-\frac{k_1}{k_2}||V||_{L^p}}\right)
\end{align*}
Now, we use that $\ua_{n-1}$ and $\ua_{n-2}$ are vector-valued and conclude
$$ ||\,|\ua_{n-1}|-|\ua_{n-2}|\,||_{L^{2p}} \leq ||\,|\ua_{n-1}-\ua_{n-2}|\, ||_{L^{2p}}
\leq ||\,\,|\ua_{n-1}-\ua_{n-2}|^2||_{L^p}^{1/2} \leq C ||\ua_{n-1}-\ua_{n-2}||_{W_p^1}. $$ Then the contractivity
constant can be bound from above by
$$L = 2k_1C^2\left(\frac{1}{2k_2}-W\right)= 2k_2\left(\frac{1}{2k_2}-W\right)=1-2k_2W =
1-\sqrt{1-4k_1k_2||V||_{L^p}}. $$
\end{proof}
\begin{remark}
Actually, it is no problem to find a suitable $\ua_0\in W_p^1(G)$ and ${\rm tr\,}\ua_0\in {\rm im\,}Q\cap
W_p^{1-1/p}(\Gamma).$ Obviously, $\ua = \vec{0}$ fulfils all necessary conditions. But we have also other
possibilities. Choose an arbitrary scalar-valued function $b(x)$ such that $||b||_{L^p}$ is small enough then
$\ua_0 = Tb$ is obviously vector-valued and fulfils all necessary conditions.
\end{remark}
\begin{remark}
More serious problem is that we had assumed the potential to be of small $L^p$-norm which need not to be fulfilled.
In general the equation
$$ \ua = T(|\ua|^2) + TV $$
is a system of nonlinear Fredholm equations of second with a weakly singular integral operator $T.$ Hence methods
to solve such types of nonlinear Fredholm equations could be applied.
\end{remark}
In what follows we will see that the solution of \ref{eqnl} is closely related to the solution of the linear
Schr\"{o}dinger equation.

The following considerations are an easy generalization from the three dimensional case which was treated in
\cite{kkw} to the n-dimensional one. We go back to the Miura transform which could be considered as a Riccati
equation as done in \cite{kkw}.\\
The Riccati equation
$$ \partial u = pu^2+qu+r, $$
where $p,\,q,\,r$ are functions, has received a lot of attention especially because of the wide range of
applications in which it appears. This equation can be reduced into its canonical form
$$ \partial w + w^2 = -v .$$
Also the Riccati equation is related to the one-dimensional Schr\"{o}dinger equation
\begin{eqnarray} \label{seq}
 -\Delta u - vu = 0
\end{eqnarray}
with a function $v$ by the substitution
$$ w = \frac{\partial u}{u} = \partial \ln u.$$
Therefore the logarithmic derivative (in the sense of Marchenko) of a scalar-valued function $u$ such that $u\not =
0$ in $G,$ is defined as
$$ \breve{\partial}u = u^{-1}Du = D(\ln u). $$
\begin{proposition}[(cf. \cite{kkw})] The scalar-valued function $\phi$ is a solution of the Schr\"{o}dinger equation
(\ref{seq}) if and only if $\ua := \breve{\partial}\phi $ is a solution of the Miura transform
$$ D\ua = v + |\ua|^2.$$
\end{proposition}
\begin{proof} Suppose that there exists a function $\phi $ such that $\ua = \breve{\partial}\phi .$ Then
$$ D\ua = \frac{1}{\phi ^2}|D\phi|^2 - \frac{1}{\phi}\Delta \phi, \quad\mbox{and}\quad
  |\ua|^2 = \frac{1}{\phi ^2} |D\phi|^2, $$
and thus
$$ D\ua - |\ua|^2 = - \frac{1}{\phi}\Delta \phi = v \iff -\Delta\phi -v\phi = 0.$$
Conversely, if $\phi $ is a solution of the Schr\"{o}dinger equation $\ua = \breve{\partial}\phi $ is a solution of the
Miura transform.
\end{proof}

\section{Application to the Gross-Pitaevskii equation}
Our considerations of the physical background of the Gross-Pitaevskii equation is based on the paper \cite{shk} and
\cite{sig}. The Gross-Pitaevskii equation for the order parameter of the mean field theory is derived by the
Fermi's zero-range pseudo-potential \cite{gross} and \cite{pit}.
\subsection{Bose-Einstein condensate}
An ideal Bose gas is a quantum-mechanical version of a classical ideal gas. It is composed of bosons, which have an
integral value of spin, and obey Bose-Einstein statistics. The statistical mechanics of bosons were developed by
Satyendra Nath Bose for photons, and extended to massive particles by Albert Einstein who realized that an ideal
gas of bosons would form a condensate at a low enough temperature, unlike a
classical ideal gas. This condensate is known as a Bose-Einstein condensate.\\[2ex]
In dilute gases, only binary collisions are taken into account. The Hamiltonian of the interacting system is:
\begin{multline*}
\hat{H}=\int dx \hat{\Psi}^{\dagger}(x)\,\left[-\frac{\hbar^2}{2m}\Delta + V(x)\right]\,\hat{\Psi}(x) + \\
   + \frac{1}{2}\int dx\,dx^{\prime} \hat{\Psi}^{\dagger}(x) \hat{\Psi}^{\dagger}(x^{\prime})U(x-x^{\prime})
   \hat{\Psi}(x^{\prime})\hat{\Psi}(x),
\end{multline*}
where $\hat{\Psi}$ is the annihilation operator, $m$ the mass of a single boson, $\hbar$ reduced Planck constant.
Experimental realizations of Boson-Einstein condensate with dilute atomic gases has been achieved in a trapping
potential $V(x)$, whose
shape is harmonic.\\[1ex]
The time evolution of the condensate is given by the Heisenberg equation
\begin{align*}
i\hbar\frac{\partial}{\partial t}  \hat{\Psi}(x,t) & = \left[  \hat{\Psi},\,\hat{H} \right] \\
  & = \left[-\frac{\hbar^2}{2m}\Delta  + V(x)+ \int \hat{\Psi}^{\dagger}(x^{\prime},t)U(x-x^{\prime})
  \hat{\Psi}(x^{\prime},t)\right]
  \hat{\Psi}(x,t).
\end{align*}
The mean-field approximation consists in replacing $\hat{\Psi}$ by the classical field $\Psi_0$. In the integral
containing $U(x^{\prime}-x),$ this is a bad approximation for short distances. However, only collisions at low
energy are relevant, which can be described by a single parameter: the s-wave scattering length $\alpha .$ Then one
can write:
$$ U(x^{\prime}-x)=g\delta (x-x^{\prime}) \quad \mbox{with} \quad g=\frac{4\pi\,\hbar^2\,\alpha }{m}. $$
And we obtain the \emph{Gross-Pitaevskii equation}
$$ i\hbar\frac{\partial}{\partial t}  \Psi_0(x,t)= \left(-\frac{\hbar^2}{2m}\Delta + V(x)+g|\Psi_0(x,t)|^2\right)
\Psi_0(x,t), $$ which is valid if
\begin{itemize}
    \item the number of atoms in the condensate $\gg 1,$
    \item low temperature
    \item $\alpha  \ll $ average distance between atoms: $\bar{\nu}\alpha ^3 \ll 1$ ($\nu$ is the density of the gas).
\end{itemize}
Here, we have the Fermi pseudo-potential (given by the Dirac distribution) has one parameter -- the scattering
length. Thus, the simplest model which can describe nonlocal repulsive interactions involves introduction of an
effective potential with two parameters: the strength of interaction and the interaction range. The symmetry
properties and the assumption of short-range interactions is enough to single out the effective potential with two
parameters. For spherically symmetric interaction the $U$ depends only on the relative distance
$U(x-x^{\prime})=U(|x-x^{\prime}|).$ Its Fourier image is then a real function of the squared wave number $\bf{k}:$
$$ \hat{U}=\hat{U}(\bf{k}^2).$$
For the Fermi pseudo-potential is $\hat{U}=1.$ In the next order of approximation we take the \emph{Lorentzian:}
$$ \hat{U}=(1+\alpha ^2\bf{k}^2)^{-1}.$$
The Lorentzian has been used recently as an approximation of the finite-range potential. In three dimensions the
Lorentzian corresponds to the \emph{Yukawa} effective \emph{potential:}
$$ U_{3D}=\frac{1}{4\pi \alpha ^2 x}exp \left(-\frac{x}{\alpha }\right), $$
and in two dimensions the effective interaction potential is given by
$$ U_{2D} = \frac{1}{2\pi \alpha ^2}K_0\left(\frac{x}{\alpha }\right), $$
where $K_0(z)$ is the MacDonald's function. \\
In the following we will use the operator form for the effective potential:
$$ \Lambda = 1-\alpha ^2 \nabla ^2 = 1- \alpha ^2 \Delta . $$
Which allows to consider the nonlocal Gross-Pitaevskii equation as a system with an additional dependent variable
$F$ describing the nonlinear term:
\begin{align}
 & \qquad (1-\alpha ^2\nabla ^2)F = |\psi |^2, \\[0.5ex]
 & \left( i\hbar\frac{\partial}{\partial t}+ \frac{\hbar^2}{2m}\Delta - V(x)- gF\right)\Psi_0(x,t)  = 0.
\end{align}

\subsection{Stationary Gross-Pitaevskii equation}\ \\
We obtain the \emph{stationary Gross Pitaevskii equation} by setting
$$ \Psi_0(x,t) = \phi(x)e^{-i\mu t/\hbar},\quad \mbox{where $\mu$ is the chemical potential,} $$
and then
\begin{align}
 &\left(-\frac{\hbar^2}{2m}\Delta + V(x)+g|\phi(x)|^2\right)\phi(x) = \mu \phi(x), \quad \mbox{or}  \nonumber \\[1ex]
 &\hspace*{3.5cm}(1-\alpha ^2\nabla ^2)F(x) = |\phi(x) |^2, \label{gpe2}  \\[1ex]
 &\left(-\frac{\hbar^2}{2m}\Delta + gF(x)+ V(x)-\mu \right)\phi(x) =  0.   \label{gpe3}
\end{align}
Where the nonlinear equation is equivalent to a system of linear equations. Here, we see that the stationary
Gross-Pitaevskii equation fits perfect in our schema. Equation (\ref{gpe2}) is a Helmholtz equation
$$ \left(-\Delta - \frac{1}{(i\alpha)^2}\right)F(x) = \frac{1}{\alpha^2}|\phi (x)|^2 $$
with constant wave number and equation (\ref{gpe3}) is a (stationary) Schr\"{o}dinger equation
$$ \left( -\Delta - \frac{-2m}{\hbar ^2}\left(gF(x)+V(x)-\mu\right)\right)\phi (x) = 0,$$
of which $\phi $ is a solution if and only if
$$ \ua = \breve\partial \phi = \phi^{-1}D\phi = D(\ln \phi) $$
is a solution of the Miura transform
$$ D\ua  - |\ua|^2 = \frac{-2m}{\hbar ^2}\left(gF(x)+V(x)-\mu\right).$$

\subsection{Non-stationary Gross-Pitaevskii equation}\ \\
Also the non-stationary Gross-Pitaevskii equation
$$ i\hbar\frac{\partial}{\partial t}  \psi(x,t)= \left(-\frac{\hbar^2}{2m}\Delta + V(x)+g|\psi(x,t)|^2\right) $$
 maybe written as a system of linear equations:
\begin{align}
 &\hspace*{3.5cm}(1-\alpha ^2\nabla ^2)F(x,t) = |\psi(x,t) |^2, \label{gpe4}  \\[1ex]
 & \left(\left(-i\hbar\frac{\partial}{\partial t}-\frac{\hbar^2}{2m}\Delta\right)+ V(x)+gF(x,t) \right) \psi(x,t) = 0.
 \label{gpe5}
\end{align}
Equation (\ref{gpe4}) is again equivalent to a Helmholtz equation with constant wave number but the second equation
(\ref{gpe5}) has due to the part $F(x,t)$ a potential which depends on the time $t.$

\end{document}